\documentclass[12pt,thmsa]{article}
\usepackage{amssymb}

\begin{document}

\author{David Carf\`{i}}
\title{Spectral expansion of Schwartz linear operators}
\date{}
\maketitle

\begin{abstract}
In this paper we prove and apply a theorem of spectral expansion for
Schwartz linear operators which have a Schwartz linearly independent
eigenfamily. This type of spectral expansion is the analogous of the
spectral expansion for self-adjoint operators of separable Hilbert spaces,
but in the case of eigenfamilies of vectors indexed by the real Euclidean
spaces. The theorem appears formally identical to the spectral expansion in
the finite dimensional case, but for the presence of continuous
superpositions instead of finite sums. The Schwartz expansion we present is
one possible rigorous and simply manageable mathematical model for the
spectral expansions used frequently in Quantum Mechanics, since it appears
in a form extremely similar to the current formulations in Physics.
\end{abstract}

\bigskip 

\bigskip

\section{\textbf{Preliminaries}}

\bigskip

In the following we shall use (with Dieudonn\'{e}) the notation $\mathcal{L}(%
\mathcal{S}_{n}^{\prime })$ for the space of $^{\mathcal{S}}$linear
endomorphisms on the space $\mathcal{S}_{n}^{\prime }$, instead of $\mathcal{%
L}(\mathcal{S}_{n}^{\prime },\mathcal{S}_{n}^{\prime })$, it is just the
space of continuous linear endomorphisms with respect to the weak* topology
(or, equivalently, with respect to the strong* topology) on the distribution
space $\mathcal{S}_{n}^{\prime }$.

\bigskip

Let $E\;$be a vector space$\;$and let $A$ be a linear operator of $E$ into $%
E $. The set of all the eigenvectors of the operator $A$ is denoted
(following Dieudonn\'{e}) by $\mathrm{E}(A)$. The set of all the eigenvalues
of the operator $A$\ is denoted by $\mathrm{ev}(A)$. Moreover, the
eigenspace corresponding to an eigenvalue $a\in \Bbb{K}$\ is denoted
(following Dieudonn\'{e}) by $\mathrm{E}_{a}(A)$. For every eigenvector $u$
of the operator $A$, there is only one eigenvalue $a$ such that $A(u)=au$,
so that we can consider the projection 
\[
e_{A}:\mathrm{E}(A)\rightarrow \mathrm{ev}(A) 
\]
associating with every eigenvector $u$ of the operator $A$ its eigenvalue,
so that $e_{A}(u)$ is the unique scalar such that 
\[
A(u)=e_{A}(u)u. 
\]
It is clear that the set $\mathrm{E}_{a}^{\neq }(A)$, collection of all
eigenvectors of $A$ corresponding to the eigenvalue $a$, coincides with the
reciprocal image $e_{A}^{-}(a)$. So that we have constructed a fiber space $(%
\mathrm{E}(A),\mathrm{ev}(A),e_{A})$, with support $\mathrm{E}(A)$, basis $%
\mathrm{ev}(A)$ and projection $e_{A}$.

\bigskip

\textbf{Remark.} We can obtain a fiber vector space, considering the set $%
E_{A}$ of all pairs $(a,u)$ in the Cartesian product $\Bbb{K\times }E$ with $%
a$ eigenvalue of the operator $A$ and such that $A(u)=au$ (so that the
partial zero-element $(a,0)$ lies in $E_{A}$, for every eigenvalue $a$ of $A$%
) with projection 
\[
e_{A}:E_{A}\rightarrow e(A):(a,u)\mapsto a, 
\]
and basis $\mathrm{ev}(A)$. Each fiber $e_{A}^{-}(a)$ of the fiber space $%
(E_{A},\mathrm{ev}(A),e_{A})$ is isomorphic (in the obvious way) with the
eigenspace $\mathrm{E}_{a}(A)$ (by the isomorphism $I_{a}:e_{A}^{-}(a)%
\rightarrow \mathrm{E}_{a}(A)$ sending the pair $(a,u)$ into the vector $u$).

\bigskip

\section{\textbf{Spectral }$^{\mathcal{S}}$\textbf{expansions}}

\bigskip

\textbf{Definition (of eigenfamily).}\emph{\ Let }$A\in \mathcal{L}(\mathcal{%
S}_{n}^{\prime })$\emph{\ be an }$^{\mathcal{S}}$\emph{linear endomorphism
on the space }$\mathcal{S}_{n}^{\prime }$\emph{, i.e. a continuous linear
endomorphism on the space }$\mathcal{S}_{n}^{\prime }$\emph{, let }$a\in 
\mathcal{O}_{M}^{(m)}$\emph{\ be a complex smooth slowly increasing function
defined on the Euclidean space }$\Bbb{R}^{m}$\emph{\ }$\emph{and}$\emph{\
let\ }$v\in \mathcal{S}(\Bbb{R}^{m},\mathcal{S}_{n}^{\prime })$\emph{\ be a
Schwartz\ family in }$\mathcal{S}_{n}^{\prime }$\emph{. We say that the
Schwartz family }$v$\emph{\ is an }$^{\mathcal{S}}$\emph{\textbf{eigenfamily
of the operator} }$A$\emph{\ \textbf{with respect to the system of
eigenvalues} }$a$\emph{\ if, for each index }$p\in \Bbb{R}^{m}$\emph{, the
vector }$v_{p}$\emph{\ is an eigenvector of the operator }$A$\emph{\ with
respect to the eigenvalue }$a(p)$\emph{. In other terms, the family }$v$%
\emph{\ is an }$^{\mathcal{S}}$\emph{eigenfamily of the operator }$A$\emph{,
with respect to the system of eigenvalues }$a$\emph{, if, for each point
index }$p\in \Bbb{R}^{m}$\emph{, we have} 
\[
A(v_{p})=a(p)v_{p}, 
\]
\emph{so that, in terms of families, the image family }$A(v)$\emph{\ can be
written as the product }$A(v)=av$\emph{, of the family }$v$\emph{\ times the
function }$a$\emph{\ (recall that the space of Schwartz families }$\mathcal{S%
}(\Bbb{R}^{m},\mathcal{S}_{n}^{\prime })$\emph{\ is a module over the ring
of slowly increasing smooth functions }$\mathcal{O}_{M}^{(m)}$\emph{).}

\bigskip

Now we can state and prove the principal theorem on spectral Schwartz
expansion.

\bigskip

\textbf{Theorem (of }$^{\mathcal{S}}$\textbf{spectral expansion).}\emph{\
Let }$A\in \mathcal{L}(\mathcal{S}_{n}^{\prime })$\emph{\ be an }$^{\mathcal{%
S}}$\emph{linear endomorphism, let }$a\in \mathcal{O}_{M}^{(m)}$\emph{\ be a
smooth slowly increasing function defined on the Euclidean space }$\Bbb{R}%
^{m}$\emph{\ }$\emph{and}$\emph{\ let }$v\in \mathcal{S}(\Bbb{R}^{m},%
\mathcal{S}_{n}^{\prime })$\emph{\ be an }$^{\mathcal{S}}$\emph{linearly
independent Schwartz eigenfamily of the operator }$A$\emph{, with respect to
the system of eigenvalues }$a$\emph{. Then, we have the spectral }$^{%
\mathcal{S}}$\emph{expansion} 
\[
A(u)=\int_{\Bbb{R}^{m}}(a[u|v])v.
\]
\emph{for each tempered distribution }$u$\emph{\ in the }$^{\mathcal{S}}$%
\emph{linear hull }$^{\mathcal{S}}\mathrm{span}(v)$\emph{\ of the
eigenfamily }$v$\emph{, where }$[u|v]$\emph{\ is the coordinate distribution
of the vector }$u$ \emph{with respect to the Schwartz basis }$v$\emph{.}

\emph{\bigskip }

\emph{Proof. }For each distribution $u$ in the $^{\mathcal{S}}$linear hull $%
^{\mathcal{S}}\mathrm{span}(v)$ of the eigenfamily $v$, we have 
\begin{eqnarray*}
A(u) &=&A\left( \int_{\Bbb{R}^{m}}[u|v]v\right) = \\
&=&\int_{\Bbb{R}^{m}}[u|v]A\left( v\right) = \\
&=&\int_{\Bbb{R}^{m}}[u|v]\left( av\right) = \\
&=&\int_{\Bbb{R}^{m}}(a[u|v])v.
\end{eqnarray*}
In fact, the third equality holds because of definition of pointwise product
of a smooth slowly increasing function by a Schwartz family, i.e. 
\begin{eqnarray*}
A(v)_{p} &=&A(v_{p})= \\
&=&a(p)v_{p}= \\
&=&(av)(p),
\end{eqnarray*}
for every point index $p$ of the family $v$, as we already have noted; and
the fourth equality holds because of elementary properties of the product of
a smooth slowly increasing function by a linear continuous operator; indeed,
for every test function $\phi $, we have 
\begin{eqnarray*}
\left( \int_{\Bbb{R}^{m}}\left[ u|v\right] (av)\right) \left( \phi \right)
&=&[u|v]((av)^{\wedge }\left( \phi \right) )= \\
&=&\left[ u|v\right] (a\widehat{v}\left( \phi \right) )= \\
&=&(a\left[ u|v\right] )(\widehat{v}\left( \phi \right) )= \\
&=&\left( \int_{\Bbb{R}^{m}}(a\left[ u|v\right] )v\right) \left( \phi
\right) ,
\end{eqnarray*}
as we well know in the general case; this concludes the proof. $\blacksquare 
$

\bigskip

\textbf{Remind (superposition of a Schwartz family with respect to an
operator).} Recall the definition of superposition of an $^{\mathcal{S}}$%
family with respect to an operator. Let $V$ be a subspace of the
distribution space $\mathcal{S}_{n}^{\prime }$, let $A\in \mathrm{Hom}(V,%
\mathcal{S}_{m}^{\prime })$ be a linear operator and let $v\in \mathcal{S}(%
\Bbb{R}^{m},\mathcal{S}_{n}^{\prime })$ be an $^{\mathcal{S}}$family of
distributions in $\mathcal{S}_{n}^{\prime }$.\emph{\ }The\emph{\ }%
superposition of the family $v$\ with respect to the operator $A$, is the
operator defined as it follows 
\[
\int_{\Bbb{R}^{m}}Av:V\rightarrow \mathcal{S}_{n}^{\prime }:u\mapsto \int_{%
\Bbb{R}^{m}}A(u)v. 
\]

\bigskip

\textbf{Remark.} So, in the conditions of the above theorem, if $V$ is the
Schwartz linear hull of the family $v$, we can write 
\[
A_{\mid V}=\int_{\Bbb{R}^{m}}a[.|v]\;v, 
\]
saying that the (domain) restriction $A_{\mid V}$ of the operator $A$ to the 
$^{\mathcal{S}}$linear hull of the family $v$ is the superposition of the
family $v$ with respect to the product of the coordinate operator 
\[
\lbrack .|v]:V\rightarrow \mathcal{S}_{m}^{\prime }:u\mapsto [u|v]=(^{t}%
\widehat{v}_{\mid (\mathcal{S}_{m}^{\prime },V)})^{-1}(u), 
\]
of the family $v$ by the system of eigenvalues $a$.

\bigskip

\textbf{Remark (on the resolution of identity).} The above theorem
generalizes the \emph{Resolution of Identity theorem}. Indeed, every $^{%
\mathcal{S}}$basis of the space $\mathcal{S}_{n}^{\prime }$ is an $^{%
\mathcal{S}}$eigenfamily of the identity operator $(.)_{\mathcal{S}%
_{n}^{\prime }}$ of the space $\mathcal{S}_{n}^{\prime }$, with respect to
the constant unitary system of eigenvalues $1_{\Bbb{R}^{m}}$, so that we
have 
\[
(.)_{\mathcal{S}_{n}^{\prime }}=\int_{\Bbb{R}^{m}}[.|v]v, 
\]
for every $^{\mathcal{S}}$basis $v$ of the space $\mathcal{S}_{n}^{\prime }$%
. Moreover, if $j_{V}$ is the canonical injection of the Schwartz linear
hull $V$ of an $^{\mathcal{S}}$linearly independent family $v$ into $%
\mathcal{S}_{n}^{\prime }$, we have 
\[
j_{V}=\int_{\Bbb{R}^{m}}[.|v]v, 
\]
where $\Bbb{R}^{m}$ is the index set of the family $v$, since the canonical
injection $j_{V}$ is just the (domain) restriction to the hull $V$ of the
identity operator $(.)_{\mathcal{S}_{n}^{\prime }}$ of the space $\mathcal{S}%
_{n}^{\prime }$.

\bigskip

\textbf{Remark (on Schwartz diagonalizable operator).} The above theorem
holds in the particular case in which there is an $^{\mathcal{S}}$basis of
the space $\mathcal{S}_{n}^{\prime }$ formed by eigenvectors of the operator 
$A$. This case is the basic theme of the theory of Schwartz diagonalizable
operators.

\bigskip

\section{\textbf{Expansions by spectral distributions}}

\bigskip

We recall that:

\bigskip

\begin{itemize}
\item  \emph{if }$A$\emph{\ is a linear continuous endomorphism of a Hilbert
space }$H$\emph{, \textbf{a spectral measure of the operator} }$A$\emph{\ is
a linear continuous operator }$\mu $\emph{\ from the space }$C^{0}(S)$\emph{
, of continuous complex functions defined on the spectrum }$S$\emph{\ of the
operator }$A$\emph{, into the space }$\mathcal{L}(H)$\emph{, of linear
continuous endomorphisms of the Hilbert space }$H$\emph{, such that the
operator }$A$\emph{\ can be seen in the integral form } 
\[
A=\mu (j_{S})=\int_{S}j_{S}\;\mu ,
\]
\emph{where }$j_{S}$\emph{\ is the canonical injection of the spectrum }$S$%
\emph{\ of the operator }$A$\emph{\ into the complex field.}
\end{itemize}

\bigskip

We shall obtain some similar spectral expansions for Schwartz-linear
operators admitting a Schwartz eigenbasis of the space.

\bigskip

In the conditions of our spectral expansion theorem (for simplicity when $v$
is a Schwartz basis of the entire space), we have 
\[
A=\int_{\Bbb{R}^{m}}(a[.|v])v, 
\]
where $a$ is the smooth slowly increasing system of eigenvalues of the
operator $A$ corresponding to the Schwartz eigensystem $v$.

\bigskip

\subsection{\textbf{Spectral distribution of a Schwartz basis}}

\bigskip

\begin{itemize}
\item  \textbf{Position of the problem.} Let $v$ be a Schwartz basis of the
space $\mathcal{S}_{n}^{\prime }$ and consider the linear operator 
\[
\mu _{v}:\mathcal{O}_{M}^{(m)}\rightarrow \mathcal{L}(\mathcal{S}%
_{n}^{\prime })
\]
defined by 
\[
\mu _{v}(f)=\int_{\Bbb{R}^{m}}(f[.|v])v,
\]
for every function $f$ in $\mathcal{O}_{M}^{(m)}$. The operator $\mu _{v}$
will be our \emph{generalized spectral measure} (or, better, our \emph{%
spectral distribution}) capable to expand (in the sense of spectral measure)
each operator $A$ admitting $v$ as a Schwartz eigenbasis.
\end{itemize}

\bigskip

\textbf{Definition (the generalized spectral measure of a Schwartz basis).} 
\emph{We call the above operator }$\mu _{v}$\emph{\ \textbf{the generalized
spectral measure of the Schwartz basis} }$v$\emph{.}

\bigskip

Note, first of all, that this operator $\mu _{v}$ is continuous, with
respect to the standard topology of the function space $\mathcal{O}%
_{M}^{(m)} $ and the pointwise topology of the operator space $\mathcal{L}(%
\mathcal{S}_{n}^{\prime })$, since it is compositions of linear continuous
operators. Namely, it is the composition of the operator chain 
\[
\mathcal{O}_{M}^{(m)}\rightarrow \mathcal{L}(\mathcal{S}_{n}^{\prime },%
\mathcal{S}_{m}^{\prime })\rightarrow \mathcal{L}(\mathcal{S}_{n}^{\prime
}), 
\]
defined by 
\[
f\mapsto f[.|v]\mapsto \int_{\Bbb{R}^{m}}(f[.|v])v. 
\]

\bigskip

\textbf{Example (the spectral distribution of the Dirac basis).} Let us
consider the Dirac family $\delta $, we have 
\begin{eqnarray*}
\mu _{\delta }(f)(u) &=&\int_{\Bbb{R}^{m}}(f[u|\delta ])\delta = \\
&=&fu= \\
&=&M_{f}(u),
\end{eqnarray*}
for every distribution $u$, so that the spectral distribution of the Dirac
family $\delta $ is the operator $M:\mathcal{O}_{M}^{(n)}\rightarrow 
\mathcal{L}(\mathcal{S}_{n}^{\prime })$ sending each slowly increasing
function $f$ of $\mathcal{O}_{M}^{(n)}$ into its multiplication operator $%
M_{f}\in \mathcal{L}(\mathcal{S}_{n}^{\prime })$.

\bigskip

We note the following trivial but meaningful properties.

\bigskip

\textbf{Proposition.}\emph{\ The spectral distribution }$\mu _{v}$\emph{\ of
a Schwartz linearly independent family }$v$\emph{\ is an algebraic linear
immersion of the space }$\mathcal{O}_{M}^{(m)}$\emph{\ into the space }$%
\mathcal{L}(\mathcal{S}_{n}^{\prime })$\emph{.}

\emph{\bigskip }

\emph{Proof.} Indeed if $\mu _{v}(f)=\mu _{v}(g)$, then 
\[
\int_{\Bbb{R}^{m}}(f[.|v])v=\int_{\Bbb{R}^{m}}(g[.|v])v, 
\]
by the Schwartz linear independence of $v$ we deduce 
\[
f[.|v]=g[.|v], 
\]
and so, for every index $p$ of the family $v$, we have 
\[
f[v_{p}|v]=g[v_{p}|v], 
\]
that is equivalent to 
\[
f\delta _{p}=g\delta _{p}, 
\]
and so, $f(p)=g(p)$ for every $p$ in $\Bbb{R}^{m}$: we have the equality $%
f=g $. $\blacksquare $

\bigskip

\textbf{Proposition. }\emph{The operator }$\mu _{v}(f)$\emph{\ has }$v$\emph{%
\ as an eigenbasis and }$f$\emph{\ as the system of eigenvalues
corresponding to }$v$\emph{.}

\bigskip

\emph{Proof.} Indeed if $\mu _{v}(f)(v_{p})=\mu _{v}(g)$, then 
\begin{eqnarray*}
\mu _{v}(f)(v_{p}) &=&\int_{\Bbb{R}^{m}}(f[v_{p}|v])v= \\
&=&\int_{\Bbb{R}^{m}}(f\delta _{p})v= \\
&=&\int_{\Bbb{R}^{m}}f(p)\delta _{p}v= \\
&=&f(p)\int_{\Bbb{R}^{m}}\delta _{p}v= \\
&=&f(p)v_{p},
\end{eqnarray*}
for every $p$ in $\Bbb{R}^{m}$: we have so 
\[
\mu _{v}(f)(v)=fv. 
\]
as we claimed. $\blacksquare $

\bigskip

The first above property is contained in the following more complete result.

\bigskip

\textbf{Proposition (of algebra homomorphism).}\emph{\ The spectral
distribution }$\mu _{v}$\emph{\ of a Schwartz basis }$v$\emph{\ is an
injective homomorphism of the function algebra }$\mathcal{O}_{M}(\Bbb{R}^{m})
$\emph{\ into the operator algebra }$\mathcal{L}(\mathcal{S}_{n}^{\prime })$%
\emph{. In particular, we have that}
\[
\mu _{v}(1_{\Bbb{R}^{m}})=(.)_{\mathcal{S}_{n}^{\prime }}.
\]

\emph{\bigskip }

\emph{Proof.} We have only to prove that 
\[
\mu _{v}(fg)=\mu _{v}(f)\circ \mu _{v}(g), 
\]
for every pair of functions $(f,g)$. Since the above operators are linear
and continuous (that is Schwartz linear) it is sufficient to prove that the
operators $\mu _{v}(fg)$ and $\mu _{v}(f)\circ \mu _{v}(g)$ are equal over
one Schwartz basis, in particular the basis $v$ itself. This is obvious,
since 
\[
\mu _{v}(fg)(v)=(fg)v, 
\]
as we already have observed in the above property and 
\[
\mu _{v}(f)\circ \mu _{v}(g)(v_{p})=\mu _{v}(f)(g(p)v_{p})=g(p)f(p)v_{p}, 
\]
for every point index $p$, for the same reason. $\blacksquare $

\bigskip

\textbf{Remark.} By the way, we note that - by the above property - each
spectral distribution $\mu _{v}$, corresponding to a Schwartz basis, is an
operator valued character of the commutative algebra $\mathcal{O}_{M}(\Bbb{R}%
^{m})$.

\bigskip

\subsection{\textbf{Operator valued spectral distributions}}

\bigskip

Let us give the formal definition of our operator valued spectral
distributions.

\bigskip

\textbf{Definition (operator valued generalized measure).} \emph{We define 
\textbf{operator valued spectral distribution} (or \textbf{generalized
spectral measure}) \textbf{of} }$\Bbb{R}^{m}$\emph{\ \textbf{into} }$%
\mathcal{L}(\mathcal{S}_{n}^{\prime })$\emph{\ any operator}
\[
\mu :\mathcal{O}_{M}^{(m)}\rightarrow \mathcal{L}(\mathcal{S}_{n}^{\prime })
\]
\emph{which is continuous with respect to the standard topology of the
function space }$\mathcal{O}_{M}^{(m)}$\emph{\ and to the pointwise topology
of the operator space }$\mathcal{L}(\mathcal{S}_{n}^{\prime })$\emph{.}

\bigskip

\textbf{Open problem.} If $\mu $ is an operator valued spectral distribution
of $\mathcal{O}_{M}^{(m)}$ into $\mathcal{L}(\mathcal{S}_{n}^{\prime })$, is
it possible to find a Schwartz basis $v$ (indexed by $\Bbb{R}^{m}$) such
that $\mu _{v}=\mu $?

\bigskip

\subsection{\textbf{Integral with respect to spectral distributions}}

\bigskip

We should, now, only define a suitable integral associated with any such
generalized measure, but this is straightforward; following the Radon
measure convention:

\bigskip

\begin{itemize}
\item  \textbf{Definition (integral with respect to an operator valued
generalized measure).} \emph{We shall put } 
\[
\int_{\Bbb{R}^{m}}f\mu :=\mu (f),
\]
\emph{for every function }$f$\emph{\ in }$\mathcal{O}_{M}^{(m)}$\emph{\ and
we shall call the value }$\mu (f)$\emph{\ \textbf{integral of the function} }%
$f$\emph{\ \textbf{with respect to the generalized measure} }$\mu $\emph{.}
\end{itemize}

\bigskip

\textbf{Remark.} Note that - at least for the moment - the juxtaposition $%
f\mu $ is not a genuine product among the function $f$ and the operator $\mu 
$.

\bigskip

\subsection{\textbf{Spectral product of operators by }$^{\mathcal{S}}$%
\textbf{families}}

\bigskip

But we should and can go further.

\bigskip

\textbf{Definition (of spectral product among operators and Schwartz
families).}\emph{\ Let }$v$\emph{\ be a Schwartz family in }$\mathcal{S}%
_{n}^{\prime }$\emph{\ indexed by }$\Bbb{R}^{m}$ \emph{and let }$B$\emph{\
be a linear continuous operator from }$\mathcal{S}_{n}^{\prime }$\emph{\
into }$\mathcal{S}_{m}^{\prime }$\emph{. We define \textbf{spectral product
of the operator} }$B$\emph{\ \textbf{by the family} }$v$\emph{, denoted by }$%
(B.v)$\emph{, as the operator-valued distribution}
\[
(B.v):\mathcal{O}_{M}^{(m)}\rightarrow \mathcal{L}(\mathcal{S}_{n}^{\prime
}):f\mapsto \int_{\Bbb{R}^{m}}(fB)v,
\]
\emph{where }$fB$\emph{\ is the product defined by}
\[
fB:\mathcal{S}_{n}^{\prime }\rightarrow \mathcal{S}_{m}^{\prime
}:(fB)(u)=fB(u),
\]
\emph{and the superposition} 
\[
\int_{\Bbb{R}^{m}}(fB)v
\]
\emph{is the superposition of the Schwartz family }$v$ \emph{with respect to
the operator }$fB$\emph{.}

\bigskip

With the above new definition we can write, eventually: 
\[
\int_{\Bbb{R}^{m}}(f[.|v])v=\int_{\Bbb{R}^{m}}f([.|v].v), 
\]
where the left hand side is a superposition of the Schwartz family $v$ with
respect to the operator $f[.|v]$ and the right hand side is the integral of
the function $f$ with respect to the operator valued distribution $([.|v].v)$
(spectral product of the coordinate operator of $v$ by the family $v$
itself).

\bigskip

\subsection{\textbf{Expansions by generalized spectral measures}}

\bigskip

Using the spectral product, our initial Schwartz spectral expansion theorem
can be written in integral form as 
\[
A=\int_{\Bbb{R}^{m}}a([.|v].v). 
\]

\bigskip

Note moreover that the formal product notation $a([.|v].v)$, appearing in
the above integral, can be (in a standard way) viewed - always - as a
genuine authentic product. Indeed,

\bigskip

\begin{itemize}
\item  \emph{define (in a standard way) the \textbf{product of a function} }$%
g$\emph{\ in }$\mathcal{O}_{M}^{(m)}$\emph{\ \textbf{by a generalized measure%
} }$\mu :\mathcal{O}_{M}^{(m)}\rightarrow \mathcal{L}(\mathcal{S}%
_{n}^{\prime })$\emph{\ as the operator }$g.\mu $\emph{\ given by} 
\[
(g.\mu )(f)=\mu (gf),
\]
\emph{for any function }$f$\emph{\ in }$\mathcal{O}_{M}^{(m)}$\emph{.}
\end{itemize}

\bigskip

The above definition is correctly given, since the product of two slowly
increasing functions in slowly increasing too.

\bigskip

\begin{itemize}
\item  \emph{Define (in a standard way) the \textbf{integral of a
generalized measure} as the value of the measure at the unitary constant
function }$1_{\Bbb{R}^{m}}$\emph{, that is as it follows} 
\[
\int_{\Bbb{R}^{m}}\mu :=\int_{\Bbb{R}^{m}}1_{\Bbb{R}^{m}}\mu .
\]
\end{itemize}

\bigskip

Then, the product $a.([.|v].v)$ is well defined and we have 
\begin{eqnarray*}
\int_{\Bbb{R}^{m}}a.([.|v].v) &=&a.([.|v].v)(1_{m})= \\
&=&([.|v].v)(a)= \\
&=&\int_{\Bbb{R}^{m}}a([.|v].v).
\end{eqnarray*}

\bigskip

So the final version of the spectral expansion theorem is 
\[
A=\int_{\Bbb{R}^{m}}a.([.|v].v), 
\]
and in the above integral does not appear any ``formal'' product but only
genuine authentic algebraic products.

\bigskip

\subsection{\textbf{Generalized spectral measure on eigenspectra}}

\bigskip

Now we desire to see spectral measures from another point of view.

\bigskip

In the conditions of our spectral expansion theorem, consider the vector
space $\mathcal{O}_{M}^{(a)}$ of all complex functions $f:S\rightarrow \Bbb{C%
}$ defined on the eigenvalue spectrum $S=a(\Bbb{R}^{m})$ of the operator $A$
and such that the composite function $f\circ a$ is a function belonging to
the space $\mathcal{O}_{M}^{(m)}$.

\bigskip

We define, for every distribution $u$ in $\mathcal{S}_{n}^{\prime }$, the
operator 
\[
\mu _{a}(u,v):\mathcal{O}_{M}^{(a)}\rightarrow \mathcal{S}_{m}^{\prime
}:f\mapsto (f\circ a)[u|v]. 
\]

\bigskip

Note that, if $j_{S}$ is the canonical immersion of the eigenvalue spectrum $%
S$ into the complex plane, we have immediately 
\[
\mu _{a}(u,v)(j_{S})=(j_{S}\circ a)[u|v]=a[u|v],
\]
so that the value of the operator $\mu _{a}(u,v)$ at the immersion $j_{S}$
is the $a$-multiple of the coordinate distribution of $u$ in the basis $v$.
We so have immediately the following superposition expansion
\[
A(u)=\int_{\Bbb{R}^{m}}\mu _{a}(u,v)(j_{S})v,
\]
for every $u$ in $\mathcal{S}_{n}^{\prime }$. But we can go further.

\bigskip

\begin{itemize}
\item  \textbf{Position of the problem.} \emph{We desire to see this
operator }$\mu _{a}(u,v)$\emph{\ as a generalized measure (always in the
sense of linear and continuous operator) on the eigenspectrum }$S$\emph{.}
\end{itemize}

\bigskip

To reach our aim we should define a right topology on the space $\mathcal{O}%
_{M}^{(a)}$ and this can be do in a standard a natural way.

\bigskip

\textbf{Topology on }$\mathcal{O}_{M}^{(a)}$\textbf{.} Let $a$ be a slowly
increasing function belonging to the space $\mathcal{O}_{M}^{(m)}$. The
function 
\[
J_{a}:\mathcal{O}_{M}^{(a)}\rightarrow \mathcal{O}_{M}^{(m)}:f\mapsto f\circ
a,
\]
is a linear injection that we call the natural injection of $\mathcal{O}%
_{M}^{(a)}$ into $\mathcal{O}_{M}^{(m)}$, so $\mathcal{O}_{M}^{(a)}$ is
linearly isomorphic with the subspace $J_{a}(\mathcal{O}_{M}^{(a)})$ of the
space $\mathcal{O}_{M}^{(m)}$. Thus $\mathcal{O}_{M}^{(a)}$ can inherit
naturally the topology of its linearly isomorphic image $J_{a}(\mathcal{O}%
_{M}^{(a)})$ via the injection $J_{a}$; a subset $O$ of $V_{a}$ is open in
this topology if and only if its image $J_{a}(O)$ - by the injection $J_{a}$
- is open in $J_{a}(\mathcal{O}_{M}^{(a)})$. This is equivalent to say that
we consider open only those sets of $\mathcal{O}_{M}^{(a)}$ which are the
reciprocal image by the injection $J_{a}$ of a open sets of the space $%
\mathcal{O}_{M}^{(m)}$. Equivalently, we endow the space $\mathcal{O}%
_{M}^{(a)}$ with the coarsest topology making the injection $J_{a}$
continuous (the so called initial topology relative to the injection $J_{a}$%
).

\bigskip

In the above conditions the operator 
\[
\mu _{a}(u,v):\mathcal{O}_{M}^{(a)}\rightarrow \mathcal{S}_{m}^{\prime
}:f\mapsto (f\circ a)[u|v], 
\]
is the composition $M_{[u|v]}\circ J_{a}$ of the continuous linear injection 
$J_{a}:\mathcal{O}_{M}^{(a)}\rightarrow \mathcal{O}_{M}^{(m)}$ with the
continuous linear multiplicative operator 
\[
M_{[u|v]}:\mathcal{O}_{M}^{(m)}\rightarrow \mathcal{S}_{m}^{\prime
}:g\mapsto g[u|v], 
\]
consequently the operator $\mu _{a}(u,v)$ is linear and continuous and so it
is a genuine generalized measure in our following sense.

\bigskip

\subsection{\textbf{Measures on eigenspectra and their integrals}}

\bigskip

\textbf{Definition (of generalized measure on parametrized spectra).}\emph{\
Let }$a$\emph{\ be a smooth slowly increasing function belonging to }$%
\mathcal{O}_{M}^{(m)}$\emph{\ and let }$S$\emph{\ be its image. We define 
\textbf{generalized measure on the set} }$S$\emph{, or more precisely 
\textbf{on the parametrization} }$a$\emph{\ \textbf{of the set} }$S$\emph{,
every linear continuous operator defined on the space }$\mathcal{O}%
_{M}^{(a)} $\emph{\ endowed with its natural topology.}

\bigskip

\textbf{Definition (integral of a generalized measure).}\emph{\ For such
generalized measures on }$S$\emph{, say }$\mu :\mathcal{O}%
_{M}^{(a)}\rightarrow \mathcal{S}_{m}^{\prime }$\emph{, we define the
integral of }$\mu $ \emph{over }$S$\emph{\ the distribution }
\[
\int_{S}\mu =\mu (1_{S}),
\]
\emph{where }$1_{S}$\emph{\ is the constant unitary function on }$S$\emph{.
Moreover, for each function }$f$\emph{\ in the algebra }$\mathcal{O}%
_{M}^{(a)}$\emph{\ we define \textbf{integral of} }$f$\emph{\ \textbf{with
respect to the generalized measure} }$\mu $\emph{\ the value }$\mu (f)$\emph{%
\ of }$\mu $\emph{\ at }$f$\emph{.}

\bigskip

The definition is well given, indeed every constant function defined on the
spectrum $S$ lives in $\mathcal{O}_{M}^{(a)}$, so it is, in particular, for
the unitary constant function $1_{S}$.

\bigskip

Moreover, we give the following.

\bigskip

\textbf{Definition (product of a generalized measure on }$a$\textbf{\ by
functions in }$\mathcal{O}_{M}^{(a)}$\textbf{).} \emph{We define product of
any function }$f$\emph{\ in the space }$\mathcal{O}_{M}^{(a)}$\emph{\ by a
generalized measure }$\mu :\mathcal{O}_{M}^{(a)}\rightarrow \mathcal{S}%
_{m}^{\prime }$\emph{\ as the operator }$f\mu $\emph{, from }$\mathcal{O}%
_{M}^{(a)}$ \emph{into }$\mathcal{S}_{m}^{\prime }$\emph{, defined by} 
\[
(f\mu )(g)=\mu (fg), 
\]
\emph{for every }$g$\emph{\ in }$\mathcal{O}_{M}^{(a)}$\emph{.}

\bigskip

\textbf{Remark.} Note that - as we have already said - the space $\mathcal{O}%
_{M}^{(a)}$ is an algebra with respect to the pointwise standard operations.
Indeed, if $f$ and $g$ lie in the subspace $\mathcal{O}_{M}^{(a)}$, then 
\begin{eqnarray*}
((fg)\circ a)(p) &=&(fg)(a(p))= \\
&=&f(a(p))g(a(p))= \\
&=&(f\circ a)(p)(g\circ a)(p)= \\
&=&(f\circ a)(g\circ a)(p),
\end{eqnarray*}
and the product of the two $\mathcal{O}_{M}^{(m)}$ functions $f\circ a$ and $%
g\circ a$ lies yet in the space $\mathcal{O}_{M}^{(m)}$.

\bigskip

\begin{itemize}
\item  \emph{The value of a generalized measure }$\mu $\emph{\ at a function 
}$f$\emph{\ in }$\mathcal{O}_{M}^{(a)}$\emph{\ (that is the integral of a
function }$f$\emph{\ in }$\mathcal{O}_{M}^{(a)}$\ \emph{with respect to the
measure }$\mu $\emph{) can be viewed as the integral of the measure }$f\mu $%
\emph{.}
\end{itemize}

\bigskip

Indeed, we have 
\[
\int_{S}f\mu =(f\mu )(1_{S})=\mu (f), 
\]
for every $f$ in $\mathcal{O}_{M}^{(a)}$.

\bigskip

\subsection{\textbf{Expansions by integration on eigenspectra}}

\bigskip

Concluding, we obtain another meaningful form of the spectral expansion
theorem.

\bigskip

\textbf{Theorem.} \emph{In the conditions of our spectral expansion theorem,
we have that the coordinate system of the image }$A(u)$\emph{\ of any
tempered distribution with respect to the Schwartz eigenbasis }$v$ \emph{of
the operator }$A$\emph{\ itself can be expanded as an integral of
generalized spectral measure; specifically we have } 
\[
\lbrack A(u)|v]=\int_{S}j_{S}\mu _{a}(u,v),
\]
\emph{for every tempered distribution }$u$\emph{.}

\emph{\bigskip }

\emph{Proof.} Indeed, we have 
\begin{eqnarray*}
\lbrack A(u)|v] &=&[\int_{\Bbb{R}^{m}}(a[u|v])v|v]= \\
&=&a[u|v]= \\
&=&\mu _{a}(u,v)(j_{S})= \\
&=&\int_{S}j_{S}\mu _{a}(u,v),
\end{eqnarray*}
for every tempered distribution $u$ in $\mathcal{S}_{n}^{\prime }$. $%
\blacksquare $

\bigskip

\subsection{\textbf{Operator valued distributions on eigenspectra}}

\bigskip

Consider now the operator 
\[
\mu _{a}[.,v]:\mathcal{O}_{M}^{(a)}\rightarrow \mathcal{L}(\mathcal{S}%
_{n}^{\prime },\mathcal{S}_{m}^{\prime }):f\mapsto (f\circ a)[.|v].
\]
This new operator can be considered also as a generalized measure on the
parametric spectrum $a$, so we have 
\[
\int_{S}j_{S}\mu _{a}[.,v]=(j_{S}\circ a)[.|v]=a[.|v],
\]
and 
\[
\int_{S}\mu _{a}[.,v]=(1_{S}\circ a)[.|v]=[.|v].
\]

\bigskip

Consider moreover the operator 
\[
\mu _{(a,v)}:\mathcal{O}_{M}^{(a)}\rightarrow \mathcal{L}(\mathcal{S}%
_{n}^{\prime }):f\mapsto \int_{\Bbb{R}^{m}}(f\circ a)[.|v]v.
\]
Also this operator can be considered as a generalized measure on the
parametric spectrum $a$, so we have 
\begin{eqnarray*}
\int_{a}j_{S}\mu _{(a,v)} &=&\int_{\Bbb{R}^{m}}(j_{S}\circ a)[.|v]v= \\
&=&\int_{\Bbb{R}^{m}}a[.|v]v= \\
&=&A,
\end{eqnarray*}
and 
\begin{eqnarray*}
\int_{a}\mu _{(a,v)} &=&\int_{\Bbb{R}^{m}}(1_{S}\circ a)[.|v]v= \\
&=&\int_{\Bbb{R}^{m}}[.|v]v= \\
&=&(.)_{\mathcal{S}_{n}^{\prime }}.
\end{eqnarray*}

Exactly the analogous of the definition of spectral measure we gave in the
beginning of the section where the space $\mathcal{O}_{M}^{(a)}$ is instead
of the space $C^{0}(S)$.

\bigskip 

\textbf{Definition (of spectral distribution of an operator).} \emph{Each
operator of the above form }$\mu _{(a,v)}$\emph{, where the pair }$(a,v)$%
\emph{\ is a Schwartz eigensolution of the operator }$A$\emph{\ (this means
simply that }$v$\emph{\ is a Schwartz eigenbasis of }$A$\emph{\ and }$a$%
\emph{\ is the corresponding system of eigenvalues) is called a \textbf{%
spectral distribution of the operator} }$A$\emph{.}

\bigskip 

\section{$^{\mathcal{S}}$\textbf{Expansions and }$^{\mathcal{S}}$\textbf{%
linear equations}}

\bigskip

Let $A$ be an $^{\mathcal{S}}$linear operator on the space $\mathcal{S}
_{n}^{\prime }$ and let $v$ be an $^{\mathcal{S}}$basis of the space $%
\mathcal{S}_{n}^{\prime }$ such that $Av=av$, with $a$ function of class $%
\mathcal{O}_{M}$. We desire to solve the $^{\mathcal{S}}$linear equation 
\[
E:A(.)=d, 
\]
with $d$ in $\mathcal{S}_{n}^{\prime }$.

\bigskip

\textbf{Theorem.}\emph{\ Let }$A$\emph{\ be an }$^{\mathcal{S}}$\emph{linear
operator on the space }$\mathcal{S}_{n}^{\prime }$\emph{\ and let }$v$\emph{%
\ be an }$^{\mathcal{S}}$\emph{basis of the space }$\mathcal{S}_{n}^{\prime }
$\emph{, indexed by the }$m$\emph{-dimensional Euclidean space, such that }$%
Av=av$\emph{, with }$a$\emph{\ function of class }$\mathcal{O}_{M}$\emph{.
Then, the }$^{\mathcal{S}}$\emph{linear equation}
\[
E:A(.)=d,
\]
\emph{with }$d$\emph{\ in }$\mathcal{S}_{n}^{\prime }$\emph{, admits (at
least) one solution if and only if the representation }$d_{v}$\emph{, of the
datum }$d$\emph{\ in the }$^{\mathcal{S}}$\emph{basis }$v$\emph{, is
divisible by the function }$a$\emph{. In this case, a solution of the
equation }$E$\emph{\ is the representation of any quotient }$q$\emph{, of
the division of }$d_{v}$\emph{\ by }$a$\emph{, in the inverse basis of }$v$%
\emph{, that is the superposition } 
\[
\int_{\Bbb{R}^{m}}qv\emph{.}
\]

\emph{\bigskip }

\emph{Proof.} $(\Rightarrow )$ Let $u$ be a solution of the equation $E$. We
have 
\[
A(u)=\int_{\Bbb{R}^{m}}a[u|v]v, 
\]
by the spectral $^{\mathcal{S}}$expansion theorem and 
\[
d=\int_{\Bbb{R}^{m}}[d|v]v, 
\]
by the definition of representation of $d$ in the basis $v$. Since $v$ is $^{%
\mathcal{S}}$linear independent, we obtain the eigen-representation of the
equality $E(u):A(u)=d$, that is the equality 
\[
a[u|v]=[d|v], 
\]
so that, the distribution $d_{v}=[d|v]$ is divisible by the function $a$,
since there exists a distribution $q$ such that 
\[
aq=d_{v}. 
\]
$(\Leftarrow )$ Vice versa, if the representation $d_{v}$ is divisible by
the function $a$, then any quotient $q$ of the division of $d_{v}$ by $a$ is
a solution of $E$. Indeed, let $q$ such a quotient, we have 
\begin{eqnarray*}
A\left( \int_{\Bbb{R}^{m}}qv\right) &=&\int_{\Bbb{R}^{m}}qA(v)= \\
&=&\int_{\Bbb{R}^{m}}q(av)= \\
&=&\int_{\Bbb{R}^{m}}(aq)v= \\
&=&\int_{\Bbb{R}^{m}}d_{v}v= \\
&=&d,
\end{eqnarray*}
as we claimed. $\blacksquare $

\bigskip

We can see an interesting application.

\bigskip

\textbf{Application (the Malgrange-Ehrenpreis theorem).} We obtain, as a
very particular case the Malgrange theorem, using the H\"{o}rmander division
of a distribution by polynomials. First of all consider that the partial
derivative $\partial _{i}$ has the Fourier basis as an $^{\mathcal{S}}$
eigenfamily, indeed we have 
\[
\partial _{i}(e^{-i(p|.)})=-ip_{i}e^{-i(p|.)}, 
\]
for every positive integer $i$ less than $n$. Consequently we have 
\[
\partial ^{j}(e^{-i(p|.)})=(-i)^{|j|}p^{j}e^{-i(p|.)}, 
\]
for every multi-index $j$; thus a differential operator $D$ with constant
coefficients, say 
\[
D=\Sigma c_{j}\partial ^{j}, 
\]
has the Fourier basis $v=(e^{-i(p|.)})_{p\in \Bbb{R}^{m}}$ as an $^{\mathcal{%
\ S}}$eigenbasis. If $q$ is the quotient of the division of a distribution $%
d $ by a polynomial $\Sigma (-i)^{|j|}c_{j}(.)^{p}$, the $^{\mathcal{S}}$%
linear equation 
\[
Du=q 
\]
has the solution 
\[
\int_{\Bbb{R}^{m}}qv, 
\]
by the above theorem, and this is exactly what the Malgrange theorem says.

\bigskip

\section{\textbf{Existence of Schwartz Green families}}

\bigskip

\textbf{Theorem.}\emph{\ Let }$L\in \mathcal{L}(\mathcal{S}_{n}^{\prime })$%
\emph{\ be an }$^{\mathcal{S}}$\emph{linear operator. Let }$\lambda $\emph{\
be an }$^{\mathcal{S}}$\emph{eigenfamily of the operator }$L$\emph{\ with
corresponding eigenvalue system }$l$\emph{, i.e. let the equality} 
\[
L(\lambda _{p})=l(p)\lambda _{p}, 
\]
\emph{hold true, for every }$p$\emph{\ in the index set, say }$I$\emph{, of
the family }$\lambda $\emph{. Assume that}

\begin{itemize}
\item  \emph{there is another }$^{\mathcal{S}}$\emph{family }$\mu $\emph{\
such that the Dirac family of the space }$\mathcal{S}_{n}^{\prime }$\emph{\
can be factorized as the product} 
\[
\mu .\lambda =\delta ,
\]
\emph{in other terms assume that the Schwartz family }$\lambda $\emph{\ has
an }$^{\mathcal{S}}$\emph{left inverse with respect to the product of
Schwartz families;}

\item  \emph{the function }$l$\emph{\ is an }$^{\mathcal{O}_{M}}$\emph{\
function, it is nowhere zero and its inverse }$l^{-1}$\emph{\ is of class }$%
\mathcal{O}_{M}$\emph{\ too, that is we require that }$l$\emph{\ is an
invertible element of the ring }$\mathcal{O}_{M}$\emph{.}
\end{itemize}

\emph{Then, the operator }$L$\emph{\ has an }$^{\mathcal{S}}$\emph{Green
family, namely the family }$G$\emph{\ defined by} 
\[
G_{p}=\int_{\Bbb{R}^{n}}(1/l)\mu _{p}\lambda , 
\]
\emph{for every index }$p$\emph{\ in }$I$\emph{.}

\emph{\bigskip }

\emph{Proof.} Indeed, for every index $p$, we have 
\begin{eqnarray*}
L(G_{p}) &=&L\left( \int_{\Bbb{R}^{n}}(1/l)\mu _{p}\lambda \right) = \\
&=&\int_{\Bbb{R}^{n}}(1/l)\mu _{p}L(\lambda )= \\
&=&\int_{\Bbb{R}^{n}}(1/l)\mu _{p}(l\lambda )= \\
&=&\int_{\Bbb{R}^{n}}l(1/l)\mu _{p}\;\lambda = \\
&=&\int_{\Bbb{R}^{n}}\mu _{p}\;\lambda = \\
&=&\delta _{p},
\end{eqnarray*}
as we claimed. $\blacksquare $

$\bigskip $

The above assumptions imply that the family $\lambda $ is a system of $^{%
\mathcal{S}}$generators for $\mathcal{S}_{n}^{\prime }$ and that the family $%
\mu $ is $^{\mathcal{S}}$linearly independent. In the particular case in
which the power $\mu .\mu $ is the factorization of the Dirac basis we
deduce that the Schwartz family $\mu $ must be a Schwartz basis too.

\bigskip

Let us generalize the preceding result.

\bigskip

\textbf{Theorem.}\emph{\ Let }$L\in \mathcal{L}(\mathcal{S}_{n}^{\prime })$%
\emph{\ be an }$^{\mathcal{S}}$\emph{linear operator with an }$^{\mathcal{S}
} $\emph{eigenfamily }$\lambda $\emph{\ and corresponding eigenvalue system }%
$l $\emph{. Assume that}

\begin{itemize}
\item  \emph{there is another }$^{\mathcal{S}}$\emph{family }$\mu $\emph{\
such that} 
\[
\mu .\lambda =\delta ,
\]

\item  \emph{any member of the family }$\mu $\emph{\ is divisible by the
function }$l$\emph{, that is there is a family }$\nu $\emph{\ of
distributions such that} 
\[
l\nu _{p}=\mu _{p},
\]
\emph{for every index }$p$\emph{\ of }$I$\emph{\ (}$\nu _{p}$\emph{\ is the
quotient of the division of }$\mu _{p}$\emph{\ by }$l$\emph{).}
\end{itemize}

\emph{Then,}

\begin{itemize}
\item  \emph{the operator }$L$\emph{\ has a Green family, namely the family }%
$G$\emph{\ defined by} 
\[
G_{p}=\int_{\Bbb{R}^{n}}\nu _{p}\lambda ,
\]
\emph{for every index }$p$\emph{.}

\item  \emph{If, moreover, the family }$\nu $\emph{\ is of class }$\mathcal{S%
}$\emph{, the operator }$L$\emph{\ has an }$^{\mathcal{S}}$\emph{Green
family, namely the family defined by the product of }$^{\mathcal{S}}$\emph{\
families }$G=\nu .\lambda $\emph{.}
\end{itemize}

\emph{\bigskip }

\emph{Proof.} Indeed, for every index $p$, we have 
\begin{eqnarray*}
L(G_{p}) &=&L\left( \int_{\Bbb{R}^{n}}\nu _{p}\lambda \right) = \\
&=&\int_{\Bbb{R}^{n}}\nu _{p}L\lambda = \\
&=&\int_{\Bbb{R}^{n}}\nu _{p}(l\lambda )= \\
&=&\int_{\Bbb{R}^{n}}(l\nu _{p})\lambda = \\
&=&\int_{\Bbb{R}^{n}}\mu _{p}\lambda = \\
&=&\delta _{p},
\end{eqnarray*}
as we claimed. $\blacksquare $

\end{document}